\documentclass[11pt]{article}
\usepackage{latexsym,graphicx}
\usepackage[utf8]{inputenc}

\usepackage{tikz,ifthen}
\usetikzlibrary{calc}

\usepackage[english]{babel}
\usepackage{amsmath}
\usepackage{amssymb,amsfonts,amsthm}
\usepackage{enumerate}
\usepackage{enumitem}
\usepackage{bbm}
\usepackage{lineno}

\usepackage{xcolor}
\usepackage{comment}

\usepackage[colorinlistoftodos]{todonotes}

\usepackage{hyperref}
\usepackage{cleveref}
\hypersetup{
	colorlinks   = true, 
	urlcolor     = black, 
	linkcolor    = black, 
	citecolor   = black 
}

\newtheorem{theorem}{Theorem}[section]

\newtheorem{lemma}[theorem]{Lemma}
\newtheorem{proposition}[theorem]{Proposition}

\newtheorem{problem}[theorem]{Problem}
\newtheorem{question}[theorem]{Question}

\newtheorem{const}[theorem]{Construction}
\newtheorem{remark}[theorem]{Remark}

\newcommand{\bA}{b^A_g}
\newcommand{\bI}{b^I_g}
\newcommand{\cA}{c^A_g}
\newcommand{\cI}{c^I_g}
\newcommand{\dang}{\mathrm{Dang}}
\newcommand{\comp}[1]{#1^C}
\newcommand{\states}{\mathcal{S}}
\newcommand{\N}{\mathbb{N}}
\newcommand{\E}{\mathbb{E}}
\newcommand{\prob}{\mathbb{P}}

\newcommand{\Mnk}{\mathcal{M}_{n,k}}
\newcommand{\barG}{\overline{G}}

\newcommand{\1}{\vspace{0.1cm}}

\colorlet{green2}{green!85!black}

\begin{document}
\title{An exploration of the balance game}

\author{$^1$Paul Dorbec,  $^2$Michael A. Henning, $^{3,4}$Zsolt Tuza, and $^5$Leo Versteegen \\
\\
$^1$Normandie Univ, UNICAEN, ENSICAEN, CNRS, GREYC,
\\ 14000 Caen, France \\
\small \tt Email: paul.dorbec@unicaen.fr \\
\\
$^2$Department of Mathematics and Applied Mathematics\\
University of Johannesburg \\
Auckland Park, 2006 South Africa \\
\small \tt Email: mahenning@uj.ac.za\\
\\
$^3$HUN-REN Alfr\'ed R\'enyi Institute of Mathematics \\
H-1053 Budapest, Re\'altanoda u.\ 13-15, Hungary \\
\\
$^4$Department of Computer Science and Systems Technology\\
University of Pannonia, Veszpr\'{e}m, Hungary \\
\small \tt Email: tuza.zsolt@mik.uni-pannon.hu \\
\\
$^5$Department of Mathematics\\
University College London \\
WC1E 6BT London, United Kingdom \\
\small \tt Email: lversteegen.math@gmail.com}

\date{}
\maketitle


\maketitle
\begin{abstract}
The balance game is played on a graph $G$ by two players, Admirable (A) and Impish (I), who take turns selecting unlabeled vertices of $G$. Admirable labels the selected vertices by~$0$ and Impish by~$1$, and the resulting label on any edge is the sum modulo~$2$ of the labels of the vertices incident to that edge. Let $e_0$ and $e_1$ denote the number of edges labeled by~$0$ and~$1$ after all the vertices are labeled. The discrepancy in the balance game is defined as $d = e_1 - e_0$. The two players have opposite goals: Admirable attempts to minimize the discrepancy~$d$ while Impish attempts to maximize~$d$. When (A) makes the first move in the game, the (A)-start game balance number, $\bA(G)$, is the value of~$d$ when both players play optimally, and when (I) makes the first move in the game, the (I)-start game balance number, $\bI(G)$, is the value of~$d$ when both players play optimally. Among other results, we show that if $G$ has order~$n$, then $-\log_2(n) \le \bA(G) \le \frac{n}{2}$ if $n$ is even and $0 \le \bA(G) \le \frac{n}{2} + \log_2(n)$ if $n$ is odd. Moreover we show that $\bA(G) + \bI(\barG) = \lfloor n/2 \rfloor$. 
\end{abstract}


\section{Introduction}

In this paper, we continue the study of the balance game played on a graph $G$ introduced and first studied by 
Krop, Mittal, and Wigal~\cite{original}. The \emph{balance game} is played by two players, Admirable and Impish, who take turns selecting the unlabeled vertices of $G$. Admirable labels the selected vertices by~$0$ and Impish by~$1$, and the resulting label on any edge is the sum modulo~$2$ of the labels of the vertices incident to that edge. Let $e_0$ and $e_1$ denote the number of edges labeled by~$0$ and~$1$ after all the vertices are labeled.  The discrepancy in the balance game is defined as $d = e_1 - e_0$. The two players have opposite goals: Admirable attempts to minimize the discrepancy~$d$ while Impish attempts to maximize~$d$. For notational convenience, sometimes we refer to Admirable as (A) and to Impish as (I). When (A) makes the first move in the game, the (A)-start game balance number, $\bA(G)$, is the value of~$d$ when both players play optimally, and when (I) makes the first move in the game, the (I)-start game balance number, $\bI(G)$, is the value of~$d$ when both players play optimally.

The balance game belongs to the growing family of \emph{competitive optimization graph games.} One of the first and best-known competitive optimization parameters is the \emph{Maker}-\emph{Breaker game} introduced in 1973 by Erd\H{o}s and Selfridge~\cite{erdos-1973}. The celebrated \emph{game chromatic number} was first introduced by Gardner~\cite{Gar-81}, and has seen extensive study, see for example~\cite{DiZh-99,KK-09,KT-94,Zhu-08}. Research has been done on competitive optimization variants of list-colouring~\cite{BorSidZs-07,Sch-09,Zhu-08}, matching~\cite{CKOW-13}, Ramsey theory~\cite{BGKMSW-11,GHK-04,GKP-08} and the domination game~\cite{BrKlRa10,BrHeKlRa-21}, among others. 
We refer the reader to the extended bibliography~\cite{Fr-12} on competitive optimization graph games, and to \cite{ZsZhu-15} for a survey on graph coloring games.

The balance game is related to the cordiality game studied in~\cite{original}. As in the balance game, the \emph{cordiality game} is played on a graph $G$ by two players, Admirable and Impish, who take turns selecting the unlabeled vertices of $G$. The game is played in a fashion similar to the balanced game, but with a modified discrepancy~$d$. The \emph{discrepancy} in the cordiality game is defined as $d = |e_1 - e_0|$.  Admirable attempts to minimize $d$ and Impish attempts to maximize~$d$. When (A) makes the first move in the game, the (A)-\emph{start game cordiality number}, $\cA(G)$, is the value of~$d$ when both players play optimally, and when (I) makes the first move in the game, the (I)-\emph{start game cordiality number}, $\cI(G)$, is the value of~$d$ when both players play optimally. The authors in~\cite{original} pointed out that $\bA(G) \le \cA(G)$ and posed the following question. 

\begin{question}
\label{Quest:Q1}
Is it true that $\bA(G) \ge 0$ for every graph $G$? 
\end{question}

We give a twofold answer to this question.
On one hand, $\bA(G)$ is non-negative for all graphs $G$ of odd order (where the \emph{order} of a graph is the number of vertices in the graph); and also $\bI(G)$ is non-negative for all graphs $G$ of even order.
On the other hand, parity really matters. For example, if $G$ is the famous Petersen graph, then $\bA(G) = -1$. 
Nevertheless, some classes of graphs admit the general lower bounds, $\bA,\bI\geq 0$ without parity constraints. An illustrative example is the class of all trees.

Our immediate aim in this paper is to obtain lower and upper bounds on the (A)-start game balance number of a general graph in terms of its order. Among other results, we show that if $G$ has order~$n$, then $-\log_2(n) \le \bA(G) \le \frac{n}{2}$ if $n$ is even and $0 \le \bA(G) \le \frac{n}{2} + \log_2(n)$ if $n$ is odd. We prove analogous results for the (I)-start game balance number, and show that if $G$ has order~$n$, then $0 \le \bI(G) \le \frac{n}{2}+\log_2(n)$ if $n$ is even and $-\log_2(n) \le \bI(G) \le \frac{n}{2}$ if $n$ is odd. We also obtain results concerning the relation between the (A)-start and (I)-start game balance numbers. For example, we show that if $G$ is a graph of order~$n$, then $\bA(G) + \bI(\barG) = \lfloor n/2 \rfloor$. Moreover, we establish lower and upper bounds on $\bI(G)$ in terms of $\bA(G)$ and the maximum degree~$\Delta$ of $G$. 
In a separate section we deal with the game balance number of paths, and observing that similar results hold also for the cordiality game, we solve in particular a further problem raised in~\cite{original}.

Already for paths, it remains unsolved to determine tight asymptotics; and more generally it is a wide open area to study the balance/cordiality game on famous graph classes, e.g.\ on planar graphs or in subclasses of perfect graphs beyond trees.

As remarked by the authors in~\cite{original}, the balance game has an interesting physical interpretation through the
well-known Ising Model from statistical physics; for example, see~\cite[Section~1.4.2]{Fr-17}  for the
relevant definitions. Two players assign spin states, one assigning ``up'', the other player ``down'', to the remaining unassigned particles in alternating turns, competing over the energy of the final configuration at game completion. As remarked in~\cite{original}, Question~\ref{Quest:Q1} through the physical interpretation may be stated as follows: Does the optimal strategy for the (A)-start balanced game keep the Hamiltonian of the model non-negative (with the assumption of an external magnetic field being absent)?

\subsection{Notation and terminology}
\label{S:notation}

For graph theory notation and terminology, we generally follow~\cite{HaHeHe-23}. Specifically, let $G$ be a graph with vertex set $V(G)$ and edge set $E(G)$, and of order $n(G) = |V(G)|$ and size $m(G) = |E(G)|$. Two vertices $u$ and $v$ of $G$ are \emph{adjacent} if $uv \in E(G)$, and are called \emph{neighbors}. The \emph{degree} of a vertex $v$ in $G$ is the number of neighbors of $v$ in $G$, and is denoted by $\deg_G(v)$, and so $\deg_G(v) = |N_G(v)|$, where $N_G(v)$ is the \emph{open neighborhood} (the set of neighbors) of $v$. The \emph{closed neighborhood} of $v$ is the set $N_G(v) = N_G(v) \cup \{v\}$. The maximum (minimum) degree among the vertices of $G$ is denoted by $\Delta(G)$ ($\delta(G)$, respectively). A \emph{cycle} on $n$ vertices is denoted by $C_n$ and a \emph{path} on $n$ vertices by $P_n$. The  \emph{complete graph} on $n$ vertices is denoted by $K_n$, while the \emph{complete bipartite graph} with one partite set of size~$n$ and the other of size~$m$ is denoted by $K_{n,m}$. For a graph $G$, let $\barG$ denote the complement of $G$, i.e., $V(\barG) = V(G)$ and two distinct vertices of $\barG$ are adjacent if and only if they are not adjacent in $G$. 

A \emph{state} of a balancing game on a graph $G$ is given by a pair $(X,Y)$ of sets of vertices of $G$, where $X$ represents the vertices selected by (A), and $Y$ represents the vertices selected by (I). In many cases we shall write $(S_0,S_1)$ for a state.

Of course, not all states are possible to arise in a game. For example, $S_0$ and $S_1$ must be disjoint, and if we play the (A)-start version, then $S_0$ may never be smaller than $S_1$. However, to avoid an excess of definitions and notations for different versions of the game, we still consider such ``impossible'' pairs of sets as states, and include them in the space of states, denoted by $\states(G)$. For a specific version of the game, we will say that a state is \emph{valid}, if it may occur in a legal game between the players, regardless of whether the state may be reached through optimal play.

The \emph{(discrepancy) score} of a state $s(S_0,S_1)$ is the number of edges between $S_0$ and $S_1$ minus the number of edges within $S_0$ and $S_1$. A \emph{strategy} is a map $\Sigma\colon \states(G)\rightarrow V(G)$ such that $\Sigma(S_0,S_1)\notin S_0\cup S_1$ unless $S_0\cup S_1 =V(G)$.

\section{General graphs}

In this section, we obtain lower and upper bounds on the (A)-start and (I)-start game balance numbers of a general graph in terms of its order. Moreover, we present results relating the (A)-start and (I)-start game balance numbers of a general graph.  

\begin{lemma}\label{lemma:symmetry}
If $G$ is a graph of order~$n$, then for any state $(S_0,S_1)$, we have $s_G(S_0,S_1)=s_G(S_1,S_0)$. Furthermore, if $S_0$ and $S_1$ partition $V(G)$, then $s_G(S_0,S_1)+s_{\barG}(S_0,S_1)=\lfloor n/2\rfloor$.
\end{lemma}
\begin{proof}
    The fact that $s_G(S_0,S_1)=s_G(S_1,S_0)$ follows from the symmetry of $S_0$ and $S_1$ in the definition of the score. For the second part, note that $s_G(S_0,S_1)+s_{\barG}(S_0,S_1)=s_{K_n}(S_0,S_1)$. It is easy to see that the latter is $\lfloor n/2\rfloor$ for any valid end state $(S_0,S_1)$.
\end{proof}

As an application of Lemma~\ref{lemma:symmetry}, we relate the (A)-start and (I)-start game balance numbers of a graph and its complement.

\begin{theorem}\label{thm:complements}
If $G$ is a graph of order~$n$, then we have $\bA(G)+\bI(\barG)=\lfloor n/2 \rfloor$.
\end{theorem}
\begin{proof}
    Let $\Sigma$ be an optimal strategy for Impish for the Admirable-start game on $G$. This means that $\Sigma$ is a function that maps a state $(S_0,S_1)$ to a vertex $v\notin S_0\cup S_1$, and if Impish plays according to this strategy, then the score at the end of the game will be at least $\bA(G)$.

    If $(S_0,S_1)$ is a valid state for the (I)-start game on $\barG$, then $(S_1,S_0)$ is a valid state in the (A)-start game on $G$. Thus, in the Impish-start game on $\barG$, Admirable may pursue the strategy
    \begin{equation*}
        \comp{\Sigma} \colon \states(G)\rightarrow V(G) \qquad (S_0,S_1) \mapsto \Sigma(S_1,S_0).
    \end{equation*}
    Suppose that after (A) plays according to $\comp{\Sigma}$, the game ends in the state $(S_0,S_1)$. By definition of $\Sigma$, we have $s_G(S_1,S_0)\ge \bA(G)$, and by \Cref{lemma:symmetry}, it follows that
    \begin{equation*}
        s_{\barG}(S_0,S_1)=s_{\barG}(S_1,S_0)=\lfloor n/2\rfloor- s_G(S_1,S_0) \le \lfloor n/2\rfloor - \bA(G).
    \end{equation*}
    Therefore, $\bI(\barG)+\bA(G)\le \lfloor n/2\rfloor$. Likewise, Impish may imitate any strategy of Admirable in the A-start game on $G$ in the I-start game on $\barG$, from which we deduce that $\bI(\barG)+\bA(G)= \lfloor n/2\rfloor$, as desired.
\end{proof}

The following proposition states that if $n$ is even, both players should prefer to start. This is because if (A) has a good strategy for the (I)-start game, they may convert this into a strategy for the (A)-start game by choosing an arbitrary vertex that they imagine (I) played in the first move. When Impish actually plays this vertex later on, Admirable may pretend that Impish played an arbitrary different vertex instead.

\begin{proposition}\label{prop:start-continuity}
If $G$ is a graph of order~$n$ with maximum degree $\Delta$, then
\[
\begin{array}{rcccll}
\bA(G) & \le & \bI(G) & \le & \bA(G)+4\Delta, & \mbox{if $n$ is even;} \1 \\
\bA(G)-2\Delta & \le & \bI(G) & \le & \bA(G)+2\Delta, & \mbox{if $n$ is odd.}
\end{array}
\]
\end{proposition}
\begin{proof}
    We begin by bounding $\bI(G)$ from above in terms of $\bA(G)$. Let $\Sigma$ be an optimal strategy for Admirable in the A-start game on $G$. Let $u=\Sigma(\empty,\empty)$ be the first move that Admirable would like to make according to this strategy and consider the following modified strategy for the I-start game on $G$.
    \begin{equation*}
        \Sigma' \colon \states(G) \rightarrow V(G) \qquad (S_0,S_1) \mapsto \begin{cases}
            \Sigma(S_0\cup \{u\}, S_1) &\text{if $u\notin S_1$},\\
            \Sigma(S_0\cup \{u\}, (S_1\setminus \{u\}) \cup \{\min \,[n]\setminus (S_0\cup S_1)\} &\text{if $u\in S_1$}.
        \end{cases}
    \end{equation*}
    If $n$ is odd, (A) may adhere to $\Sigma'$ throughout the whole game, and for the final state $(S_0,S_1)$ of the game, we must have $u\in S_1$. Because $\Sigma$ is optimal, it follows that $s_G(S_0 \cup \{u\}, S_1\setminus \{u\})\le \bA(G)$, and because changing the label of $u$ cannot change the score by more than twice the degree of $u$, we conclude that $\bI(G)\le s_G(S_0,S_1)\le \bA(G)+2\Delta$.

    If $n$ is even, (A) may adhere to $\Sigma'$ up until the final move. Let $v$ denote the last unlabeled vertex, which (A) must play in the final move. If $v=u$, then the final state of the game that (A) imagined game coincides with the final state $(S_0,S_1)$ of the actual game so that $s_G(S_0,S_1)\le \bA(G)$. If $v\neq u$, then swapping the labels of $v$ and $u$ will change the score by at most $4\Delta$, which means that $\bI(G)\le \bA(G) + 4\Delta$.

    To bound $\bI(G)$ from below, note that Impish may employ an analogous strategy modification as Admirable above, which shows that $\bI(G)\ge \bA(G)-4\Delta$ if $n$ is even, and $\bI(G)\ge \bA(G)-2\Delta$ if $n$ is odd. However, if $n$ is even, Impish can do better. Let $\Sigma$ be a strategy of Impish for the (A)-start game on $G$, and assume without loss of generality that $V(G)=[n]$. Consider the following strategy for (I) in the (I)-start game on $G$.
    \begin{equation*}
        \Sigma'\colon \states(G) \rightarrow [n] \qquad (S_0,S_1)\mapsto \Sigma(S_0\cup \{\min V(G)\setminus (S_0,S_1)\},S_1).
    \end{equation*}
    Put in words, Impish imagines that Admirable starts by playing the vertex 1. After that, when it is not the last move, and Admirable plays the vertex $v=\min V(G)\setminus(S_0,S_1)$, Impish imagines that Admirable played the vertex $\min V(G)\setminus(S_0,S_1\cup \{v\})$ instead. On the last move, Admirable has no choice but to play $\min V(G)\setminus(S_0,S_1)$, and the end state of the game coincides with the end state of the (A)-start game that Impish imagined.
\end{proof}

In general, it will also be useful to observe that the balancing game numbers are ``continuous'' in the set of edges.

\begin{proposition}
If $G$ and $H$ are two graphs on the same vertex set, such that the symmetric difference of $E(G)$ and $E(H)$ has size $k$, then $\vert \bA(G)-\bA(H)\vert\le k$ and $\vert \bI(G)-\bI(H)\vert\le k$.
\end{proposition}
\begin{proof}
    Either player may play in $H$ as in $G$ and the score at the end differs at most by the number of edges.
\end{proof}

We are now in a position to present lower and upper bounds on the (A)-start and (I)-start game balance numbers of a general graph in terms of its order~$n$. 

\begin{theorem}   \label{thm:global-bounds}
If $G$ is a graph of order~$n$, then we have
\[
\begin{array}{rcccll}
-\log_2(n) & \le & \bA(G) & \le & \frac{n}{2}, & \mbox{if $n$ is even;} \1 \\
0 & \le & \bA(G) & \le & \frac{n}{2} + \log_2(n), & \mbox{if $n$ is odd;} 
\end{array} 
\]
and 
\[
\begin{array}{rcccll}
0 & \le & \bI(G) & \le & \frac{n}{2}+\log_2(n), & \mbox{if $n$ is even;} \1 \\
-\log_2(n) & \le & \bI(G) & \le & \frac{n}{2} & \mbox{if $n$ is odd.}
\end{array}
\] 
\end{theorem}
\begin{proof}
    It is enough to show the bounds for the (A)-start game when $n$ is even and the (I)-start game when $n$ is odd as once these are established, the other bounds follow readily from \Cref{thm:complements}.
    
    We first consider the case of the (A) start game for even $n$. Recall the definition of a state $(S_0,S_1)$ of the game and its score, which is the number of edges between $S_0$ and $S_1$ minus the number of edges within $S_0$ and $S_1$. 
    
    For a vertex $v$ and a game state $(S_0,S_1)$, the \emph{score of $v$ in state $(S_0,S_1)$} is defined as $s(v;S_0,S_1)=\vert N[v] \cap S_0\vert - \vert N[v] \cap S_1\vert$, where recall that $N[v]$ is the closed neighborhood of $v$. Thus, we have that
    \begin{align*}
        s(S_0\cup \{v\},S_1)&=s(S_0,S_1)-s(v;S_0,S_1)&
        &\text{and}&
        s(S_0,S_1\cup \{v\})&=s(S_0,S_1)+s(v;S_0,S_1).
    \end{align*}

    Note that $s(S_0,S_1)$ counts precisely the number of unbalanced edges minus the number of balanced edges. Therefore, to show that $\bA(G)\le n/2$, if $n$ is even, it is enough to show that Admirable has a strategy so that $s(S_0,S_1)\le n/2$ at the end of the game. Indeed, if it is Admirable's turn to play from state $(S_0,S_1)$, they simply play an unlabeled vertex $v$ for which $s(v;S_0,S_1)$ is maximum. This moves the game into the state $(S_0\cup \{v\},S_1)$, and since $n$ is even, this was not the last move of the game and (I) has to play next. For every unlabeled $w\in V(G)$, we have
    \begin{equation*}
        s(w;S_0\cup \{v\},S_1) \le s(w;S_0,S_1)+1.
    \end{equation*}

    Thus, by choice of $v$, whichever vertex $w$ Impish plays next, we have that $s(w;S_0\cup \{v\},S_1)\le s(v;S_0,S_1)+1$, and therefore,
    \begin{equation*}
        s(S_0\cup \{v\},S_1\cup \{w\}) \le s(S_0,S_1) - s(v; S_0, S_1)+ s(w;S_0\cup \{v\},S_1) \le s(S_0,S_1)+1.
    \end{equation*}
Hence, in two moves the score increases by at most 1, and
    it follows by induction that for every even $k\le n$, after $k$ moves were played, we reach a state $(S_0,S_1)$ such that $s(S_0,S_1)\le k/2$. In particular, $s(S_0,S_1)\le n/2$ once $S_0\cup S_1 = V(G)$.

\medskip

    We now move on to the lower bound. Define the \emph{danger potential} of a state $(S_0,S_1)$ as
    \begin{equation*}
        \dang(S_0,S_1)=2^{-s(S_0,S_1)}\sum_{v\notin S_0\cup S_1} 2^{-s(v;S_0,S_1)}.
    \end{equation*}

    We show that Impish can play in a way such that for every state $(S_0,S_1)$ after an odd number of moves $k<n$, we always have $\dang(S_0,S_1)\le n-1$. Note that this implies that $\bA(G)> -\log_2 n$. Indeed, if $v$ is the only remaining unlabeled vertex in the state $(S_0,S_1)$ before the final move, we have that
    \begin{equation*}
        -s(S_0,S_1) - s(v;S_0,S_1) = \log_2 \left(\dang(S_0,S_1) \right) < \log_2 n.
    \end{equation*}

    It follows that the state $(S_0,S_1\cup \{v\})$, which we arrive at after (I) plays $v$, has score at least $-\log_2 n$.

    It remains to confirm that Impish has a strategy to ensure that $\dang(S_0,S_1)\le n-1$ whenever an odd number of vertices has been played. Note first that the inequality holds trivially after the first move of Admirable. Assume the inequalities hold in state $(S_0,S_1)$ and it is Impish's turn to play. Let $U$ be the set of vertices $u$ that maximize $s(u;S_0,S_1)$ and denote this maximum by $m$. Among the elements of $U$, Impish plays the vertex $u$ that minimizes the quantity
    \begin{equation*}
        r(u):=\sum_{w\in N(u) \setminus (S_0\cup S_1)} 2^{-s(w;S_0,S_1)},
    \end{equation*}
    and Admirable responds by playing another vertex $v$. If $v\notin U$ or if $v$ is a neighbor of $u$, then $s(v;S_0,S_1\cup \{u\}) \le s(v;S_0,S_1)<m$, and therefore, $s(S_0\cup \{v\},S_1\cup \{u\})\ge s(S_0,S_1)+1$ so that
    \begin{align*}
        \dang(S_0\cup \{v\},S_1\cup \{u\}) \le 2^{-s(S_0,S_1)-1} \sum_{\substack{w\notin S_0\cup S_1\\w\notin \{u,v\}}} 2^{-s(w;S_0\cup \{v\},S_1\cup \{u\})+1}<\dang(S_0,S_1).
    \end{align*}
    If $v\in U$ and $u$ and $v$ are not adjacent, we must have $r(v)\ge r(u)$, and hence,
    \begin{align*}
        \dang(S_0\cup \{v\},S_1\cup \{u\}) < \dang(S_0,S_1) + 2^{-s(S_0,S_1)} (r(u)-r(v))\le \dang(S_0,S_1),
    \end{align*}
    completing the proof for the (A)-start game when $n$ is even.

    The proof for the (I)-start game when $n$ is odd is very similar. For the upper bound, observe that after the first move of (I), the score is still~$0$. Thus, it is still enough for (A) to play the vertex $v$ that maximizes $s(v;S_0,S_1)$ in every move. For the lower bound, (I) plays such that $\dang(S_0,S_1)\le n$ after every \text{even} number of moves, which leads to the same analysis for the last move as in the (A)-start case for even $n$, except that we have the immaterially weaker bound $\dang(S_0,S_1)\le n$ instead of $\dang(S_0,S_1)\le n-1$.
\end{proof}

\subsection{Graphs with prescribed balancing number}\label{subsection:prescribed}

In Theorem \ref{thm:global-bounds} we proved that $\bA(G)$ cannot exceed $n/2$ for any graph $G$ of even order $n$.
Many nice graphs show that this bound is tight, for instance
\begin{itemize}
    \item complete graphs of even order;
    \item diamond necklaces;
    \item cycles of triangles (that extend each edge of $C_{n/2}$ to a triangle by an external vertex of degree 2).
\end{itemize}
Validity for these examples will be implied by the very general construction given below.

\begin{const}
For $n \ge 2k$, let $f_1,\dots,f_{k}$ be $k$ copies of $K_2$,
if $n \ge 2k+2$, let \(f_{k+1}, \ldots f_{\lfloor\frac{n}{2}\rfloor}\)
be copies of $\overline{K_2}$, and add a vertex $v$ if $n$ is odd.
Between any two $f_i,f_j$ insert either no edges, or a $P_3$, or a $K_{2,2}$.
Between $v$ and any $f_i$, also insert either no edges or a $P_3$.
Denote by $\Mnk$ the class of graphs obtained in this way.
\end{const}

\begin{theorem}   \label{t:mnk}
For every graph $G\in\Mnk$ we have $\bA(G) = \bI(G) = k$.
\end{theorem}
\begin{proof}
Let $G\in \Mnk$.
First observe that if all pairs $f_i$ in $G$ are colored with both colors~$0$ and~$1$, then all edges between two $f_i,f_j$ or between $v$ and some $f_i$
have scores summing up to~$0$. On the other hand, all the edges within the $k$ pairs formed on $K_2$ contribute $k$, namely one for each such pair. Thus, a game ending up with such a coloring of $G$ will have score $k$.

To prove that $\bA(G) = \bI(G) = k$, we design a rather simple strategy for each player enforcing that score.

When $n$ is even, our strategy as the second player (whether it is (A) or (I)) is simple.
Whenever our opponent colors a vertex in a pair,
we color the other vertex of the pair (with the opposite color). The resulting coloring ensures a score of~$k$.

As first player, our strategy is similar.
We play any vertex as our first vertex, opening one pair.
If our opponent plays in a different (new) pair, we answer in that new pair.
If our opponent plays the second vertex of the open pair, we play any other vertex,
opening a new pair.
At the end of the game, our opponent has to conclude the game
playing the second vertex of the open pair
and complete the required coloring.

When $n$ is odd,
our strategy as the first player is to start on vertex $v$ and then adopt
our second player's strategy above. As the second player, we follow moves
of our opponent unless the vertex played is $v$,
in which case we adopt our above first player strategy. If our opponent never plays $v$, we will simply finish the game playing $v$ once all the pairs received different colors.
\end{proof}

The result of Theorem~\ref{t:mnk} proves that there are infinitely many graphs of order at least $2k$ having score $k$ when $k$ is non-negative.
In particular, when $n$ is even and $k=\frac{n}{2}$, $\Mnk$ forms a large family
achieving the upper bound of Theorem~\ref{thm:global-bounds}.
On the other hand, when $n$ is odd and $k=0$, $\Mnk$ forms a large family
achieving the lower bound of Theorem~\ref{thm:global-bounds}.

We don't have families tightening the other bounds of Theorem~\ref{thm:global-bounds}
(i.e., the lower bound for even $n$ and the upper bound for odd $n$), though we know that those bounds could not be matched to the bounds for the other parity.

In particular, it should be noted that the (A)-start game balance number of a graph of even order can be negative, which answers \emph{Open question 2} from \cite{original}.

\begin{theorem}
    For the famous Petersen graph $P$, we have
    \( \bA(P) = -1\).
\end{theorem}
\begin{proof}
The strategy of (A) is to build a connected induced subgraph $H$ of order~$5$. This can successfully be achieved because from any subtree $T$ of order $k$ there are \(3k - 2(k-1) = k+2\) edges going to $P-T$, and since $P$ has girth~$5$, those edges determine either $k+2$ or $k+1$ distinct vertices (the latter may occur if $k=4$ and the subtree is $P_4$). Hence, using only $k$ vertices, (I) cannot block all neighbors of $T$. Therefore, (A) has a choice for the connected extension of $T$, eventually yielding an $H$ as claimed.

If $H$ is $C_5$, then also the subgraph selected by (I) is a $C_5$, hence $e_1 = 5$ and $e_0 = 10$, $d = -5$. If $H$ is a tree, then $k+2 = 5+2 = 7$ edges go out from $H$, hence $e_1 = 7$ and $e_0 = 15 - 7 = 8$, $d = -1$.
\end{proof}

Observe also that the $5$-cycle is a graph of odd order having $\bA(G)$ larger than $\frac{n}{2}$. Indeed, $\bA(C_5) = 3$ (while $\bI(C_5) = -1$).

However, we were not able to disprove a slightly weaker bound, leading to the question:

\begin{problem}
    Determine whether the inequality $-1 \le \bA(G) \le \lceil \frac{n}{2}\rceil$ holds for all graphs $G$ of order~$n$.
\end{problem}

\subsection{Other general bounds}

In this section, we present other general bounds on the (A)-start and (I)-start game balance numbers. 

\begin{lemma} \label{l:imagination}
    Given a graph $G$ of maximum degree $\Delta$ and one of its vertices $v$,
    we have
    \[
    \begin{cases}
    \bI(G) \le \bA(G-v) + 2\Delta + \deg(v), \\
    \bA(G) \ge \bI(G-v) - 2\Delta - \deg(v).
    \end{cases}
    \]
\end{lemma}
\begin{proof}
    Let us first describe the strategy for Admirable when playing second on $G$. 
    The key idea is that (A) follows the strategy for playing first on $G-v$, ignoring (I)'s first move (say on some vertex $u$)
    as well as ignoring the existence of $v$.
    Obviously, (I) won't color $u$ a second time.
    What may happen is that (I) colors $v$ in $G$,
    in which case (A) imagines (I) just colored $u$ instead,
    and (A)'s imagined game will fit closely to the real game restricted to $G-v$.
    Another possibility is that instead, (A)'s strategy would be to color $u$,
    which is not a legal move in the real game.
    In that second case, (A) colors the vertex $v$ instead, and continues the game
    as if $u$ was colored.

    When the game gets to its end, the edges have the same status, except possibly for the at most $\Delta$ edges incident to $u$, and the $\deg(v)$ edges incident to $v$.
    The difference in the score may be that $\Delta$ monochromatic edges become bi-colored (resulting in a $+2\Delta$ difference in the score), and $\deg(v)$ additional bi-colored edges are added. Therefore, the score is at most \(\bA(G-v) + 2\Delta + \deg(v)\), and the first bound follows.

    The second bound can be obtained following the same idea for Impish.
\end{proof}

\section{Classical families of graphs}

\begin{lemma}\label{l:2-extra-vertices}
If $G$ is a graph with two vertices $u$ and $u'$ such that
\begin{itemize}
    \item $u$ has exactly one neighbor in $G - \{u'\}$, denoted $v$,
    \item $u'$ is only adjacent to either $u$ or $v$
\end{itemize}
then $\bA(G) \ge \bA(G - \{u,u'\})$ and $\bI(G)\ge \bI(G-\{u,u'\})$.
\end{lemma}
\begin{proof}
    Let $G' = G - \{u,u'\}$.
    For showing the desired lower bound on $\bA(G)$, we describe a strategy for Impish.
    Impish plays $G$ as if on $G'$, except that
    if Admirable plays on one of $u$ or $u'$, then Impish answers on the other,
    and the game will continue as if on $G'$.
    If the game ends on $G'$ before $u$ and $u'$ are played, then Impish may be forced to play the first move on $u$ or $u'$, and then Admirable must play
    the second. In both cases, $u$ and $u'$ are given different colors,
    though possibly $u$ gets the same color as $u'$.
    As a consequence, $e_1(G) \ge e_1(G') +1 $ and $e_0(G) \le e_0(G') + 1$, and so the resulting score is no less than $\bA(G')$, as desired.

    The proof of the bound for $\bI(G)$ follows along the same line.
\end{proof}

\begin{figure}
    \centering
    \begin{tikzpicture}
 [thick,scale=1,
	vertex/.style={circle,draw,fill=white,inner sep=0pt,minimum size=2mm},
	]

	\path (0,0) coordinate (v) ++(0:1) coordinate (u) ++(0:1) coordinate (up);
    
	\draw (v) -- (u) -- (up);
    
    \draw (u) node[vertex,label=90:$u$] {};
    \draw (up) node[vertex,label=90:$u'$] {};  
    \draw (v) node[vertex,label=90:$v$] {}; 
        
    \draw  (-1,0) node {$G'$}; 

	\draw (-1,0) circle[radius=1.5] ;

	\path (6,0) coordinate (v2) +(30:1) coordinate (u2) +(-30:1) coordinate (up2);
    
	\draw (v2) -- (u2)  (v2) -- (up2);
    
    \draw (u2) node[vertex,label=90:$u$] {};
    \draw (up2) node[vertex,label=90:$u'$] {};  
    \draw (v2) node[vertex,label=90:$v$] {}; 
        
	\draw (v2) ++ (-1,0) circle[radius=1.5] ;

	\draw  (v2) ++ (-1,0) node {$G'$}; 

\end{tikzpicture}
    \caption{Cases for Lemma~\ref{l:2-extra-vertices}}
    \label{f:2-extra-vertices}
\end{figure}
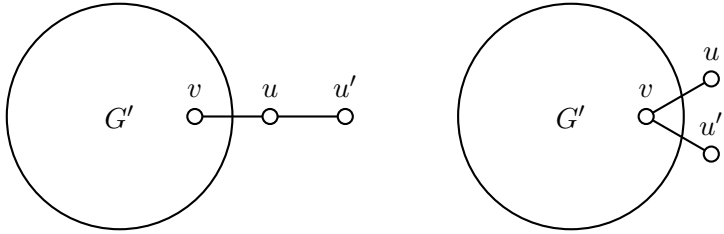

We remark that part of the proof in the above lemma was used already in~\cite{original} to prove that $\bA(P_n)$ is non-negative.

As a consequence of this lemma, we get the following result:

\begin{theorem}
    If $T$ is a tree, then $\bA(T)\ge 0$ and $\bI(T) \ge 0$.
\end{theorem}
\begin{proof}
    The proof goes by induction. The base cases are the trees of order~$1$ or~$2$, which necessarily have scores $0$ and $1$, respectively.
    Suppose $T$ is a larger tree. 
Any longest (diametric) path ends with a pendant $P_3$ or its last but one vertex is adjacent to more than one leaf.
    We recognize the two possible cases dealt with by Lemma~\ref{l:2-extra-vertices}.
\end{proof}

\begin{proposition}   \label{p:complete}
    For complete graphs $K_n$ we have
    \[\bA(K_n) = \bI(K_n) = \left\lfloor \frac{n}{2}\right\rfloor\]
\end{proposition}
\begin{proof}
The claim follows immediately from \Cref{thm:complements} and the fact that for the empty (edgeless) graph $E_n=\overline{K_n}$, we have $\bA(E_n)=\bI(E_n)=0$.
Alternatively, at the end of any game on the complete graph,
    there will be $\lfloor\frac{n}{2}\rfloor$ vertices of one color and
    $\lceil\frac{n}{2}\rceil$ vertices of the other color.
    A simple calculation brings to the conclusion.
\end{proof}

We note that the above formula is also a very particular case of Theorem \ref{t:mnk}.

\begin{proposition}   \label{p:comp-bip}
    For a complete bipartite graph $G=K_{p,q}$ on $n=p+q$ vertices,
    \[\bA(G) = \bI(G) =
    \begin{cases}
    0, & \mbox{ if $pq$ is even;}\\
    1, & \mbox{ if $pq$ is odd.}\\
    \end{cases}
    \]
\end{proposition}
\begin{proof}
    Let $G$ be a complete bipartite graph, with partite sets $V_1$ and $V_2$.
    We prove first that whatever strategies are chosen,
    the score will end up non-negative. At the end of the game, we  denote the number of vertices colored~$0$ in $V_1$ by $a$, the number of vertices colored~$1$ by $a+x$ ($x$ may be negative), and the number of vertices colored~$1$ in $V_1$ by $b$.
    Since the turns alternate starting by Admirable,
    the number of vertices colored $0$ in $V_2$ must be $b+x + \varepsilon$ where $\varepsilon=n\bmod 2$.
    We compute the score
    \[
    \begin{array}{rcl}
    e_1 - e_0 & = & ab + (a+x)(b+x+\varepsilon) - a(b+x+\varepsilon) - (a+x)b\\
     & = & x(x+\epsilon) = x^2 + \varepsilon x
     \end{array}
    \]
    Since $x$ is an integer, $x^2 \ge |x|$, and so the latest is non-negative and this concludes the statement. 

For the upper bound, Admirable will ensure $x=0$ if at least one of $p$ and $q$ is even, or $x=1$ if both are odd (but then $\varepsilon=0$).
This is achieved by the strategy that (A) plays after each move of (I) in the same vertex class, as long as possible.
    If $n$ is odd, (A) makes the first move in the odd vertex class. The computation thus concludes the proof.     
    \end{proof}

\subsection{Random graphs}

Let $G_n\sim G(n,1/2)$ be an Erdős–Rényi–Gilbert random graph. By \Cref{thm:complements}, the expectation of $(\bA(G_n)+\bI(G_n))$ (or, equivalently, the average of $\bA(G_n)$ and $\bI(G_n)$) can be determined.

\begin{theorem}
\begin{equation*}
    \E_{G_n} [\bA(G_n)+\bI(G_n)]=\E_{G_n} [\bA(G_n) + \bI(\barG_n)]= \frac{\lfloor n/2\rfloor}{2}.
\end{equation*}
\end{theorem}

If $n$ is even, \Cref{prop:start-continuity} tells us that $\bA(G_n)\le \bI(G_n)$, so that we must have 
\begin{align*}
    \E_{G_n} [\bA(G_n)]\le n/4 \le \E_{G_n} [\bI(G_n)].
\end{align*}

\begin{problem}
    Is it true that
    \begin{equation*}
        \lim_{n\rightarrow \infty}  \frac{\E_{G_{2n}}[\bA(G_{2n})]}{2n}=\lim_{n\rightarrow \infty}  \frac{\E_{G_{2n}}[\bI(G_{2n})]}{2n}=\frac{1}{4}\,?
    \end{equation*}
    Equivalently, is $\bA(G)=\bI(G)+o(1)$ with high probability if $G$ has an even number of vertices?
\end{problem}

Once the expectations are established for graphs of even order, it seems likely that one can approximate the expectations for graphs of odd order by deleting a single vertex.

It would be interesting to see how far the random variables $\bA(G)$ and $\bI(G)$ are concentrated around their mean. By \Cref{thm:global-bounds}, the variable $X=\bA(G_n)+\log_2 n$ is always positive so that we may apply Markov's inequality to get a bound on the probability that $\bA$ is large when $n$ is even, yielding

\begin{equation}\label{eq:prob-bound}
    \prob\left(\bA(G_n) \ge (1/4+c)n\right) \le \frac{n/4+\log_2 n}{(1/4+c)n} = \frac{1}{1+4c}+o(1).
\end{equation}

Note, however, that this bound is very weak: we know already from \Cref{thm:global-bounds} that $\bA(G_n)/n$ is surely bounded by $1/2$, but \eqref{eq:prob-bound} shows only that $\prob(\bA(G_n) \ge n/2)\le 1/2$. In principle, one may use (for example) an edge-exposure martingale to show somewhat stronger concentration results. However, since $\Theta(n^2)$ edges are exposed throughout the process, we cannot hope to obtain a standard deviation of order $o(n)$ in this way. Furthermore, since we have shown in \Cref{subsection:prescribed} that for all $k\in [0,n/2]$ there exists a relatively large class of graphs for which $\bA(G)=\bI(G)=k$, it seems at least plausible that $\bA(G_n)$ and $\bI(G_n)$ are in fact distributed fairly uniformly.

\begin{problem}
    Determine whether $\bA(G_n)$ and $\bI(G_n)$ are tightly concentrated around their mean.
\end{problem}

\section{Paths}

We computed the small values of the (A)-start and (I)-start game balance numbers with a computer program to obtain the results presented in Table~\ref{t:scores}.

\begin{table}[htb]\begin{center}
\begin{tabular}{|l|cccccccccccccccccc|}
 \hline
 $n$ & 1 & 2 & 3 & 4 & 5 & 6 & 7 & 8 & 9 & 10 & 11 & 12 & 13 & 14 & 15 & 16 & 17 & 18 \\
 \hline
 (A)-start & 0 & 1 & 0 & 1 & 2 & 1 & 2 & 1 & 2 & 1 & 4 & 1 & 4 & 1 & 4 & 3 & 4 & 3 \\
 (I)-start & 0 & 1 & 0 & 1 & 0 & 1 & 0 & 3 & 0 & 3 & 0 & 3 & 2 & 3 & 2 & 3 & 2 & 5 \\
\hline
\end{tabular}
\caption{scores of the game for small paths}\label{t:scores}
\end{center}
\end{table}

\begin{lemma}
For paths $P_n$ of order~$n$, we have
    \[
    \begin{cases}
    \bI(P_{n+1}) \ge \bA(P_n) - 1, \\
    \bA(P_{n+1}) \le \bI(P_n) + 1.
    \end{cases}
    \]
\end{lemma}

\begin{proof}
    Impish, playing first on $P_{n+1}$, can color the last vertex of the path and
    then follow an (A)-start strategy to continue the game on the $n$ first vertices
    of the path, imagining the last vertex does not exist. At the end, maybe the
    edge $v_nv_{n+1}$ is monochromatic, and so the resulting score is at least \(\bA(P_n) - 1\).

    As for Admirable, they can follow the corresponding strategy when starting on
    \(P_{n+1}\), relying on an optimal (I)-start strategy on \(P_n\).
    The last edge may be bi-colored and thus the score is at most \(\bI(P_n) + 1\).
\end{proof}

Note that the first inequality of this lemma could also have been inferred
from the second and Lemma~\ref{l:2-extra-vertices} which implies that $\bI(P_{n+2})\ge \bI(P_n)$.

Moreover, we note that thanks to this lemma and Lemma \ref{l:2-extra-vertices},
one can infer from Table~\ref{t:scores} that $\bA(P_{17}) = 4$ without any computation.

\begin{lemma}   \label{l:rec4}
For paths $P_n$ of order~$n$, we have
    \[
    \begin{cases}
    \bA(P_{n+4}) \le \bA(P_n) +2,\\
    \bI(P_{n+4}) \le \bI(P_n) + 2.
    \end{cases}
    \]
\end{lemma}
\begin{proof}
    Admirable plays on $P_{n+4}$ as if on $P_n$, ignoring
    the $P_4$ at some extremity of the path.
    If at some point, (I) plays a vertex $u$ on the $P_4$,
    (A) answers by playing the only vertex on the $P_4$ which
    is at distance 2 from $u$. Then, the two edges from $u$ and $v$
    to their common neighbor will compensate whatever happens later.
    Later on, if (I) colors a third vertex in the $P_4$,
    (A) colors the fourth so that the players alternation on the main
    part of the path is not altered.

    If (I) totally ignores the $P_4$, and (A) is eventually forced to play first on it,
    then (A) plays one of the middle vertices $u$ of the $P_4$,
    so whatever vertex (I) colors, then
    (A) can color another vertex on the $P_4$ adjacent to $u$.

    We infer that at the end of the game, there are at most three extra bi-colored edges and at least one extra bi-colored edge. The resulting score is thus at most~$2$ more than the score on the path on $n$ vertices.
\end{proof}

\begin{remark}   \label{r:rec6}
In a similar way, using the fact that
 $\bA(P_6) = \bI(P_6) = 1$, the inequalities
    \[
    \begin{cases}
    \bA(P_{n+6}) \le \bA(P_n) +2, \\
    \bI(P_{n+6}) \le \bI(P_n) + 2
    \end{cases}
    \]
can also be obtained by considering optimal strategies on $P_6$
at an extremity of $P_{n+6}$,
with a more complicated description of the strategy.
\end{remark}

\begin{lemma}   \label{l:rec16}
The following holds:
\begin{enumerate}[label=\rm(\alph*)]
    \item For $n$ even, $\bA(P_n) \ge 2\lfloor\frac{n}{16}\rfloor + 1$.
    \item For $n$ odd,  $\bA(P_n) \le 4\lceil\frac{n}{16}\rceil$.
    \end{enumerate}
\end{lemma}

\begin{proof}
    We first show the lower bound (a) for even paths.
    We use the fact stated in Table~\ref{t:scores} that $\bA(P_{16}) = \bI(P_{16}) = 3$, and also that for every even path $P_r$ with even $r < 16$, $\bA(P_r) = 1$.
    Impish's strategy consists in seeing the path on $n = 16k +r$ vertices as
    a collection of $k$ paths on 16 vertices and (possibly) a small remaining path
    on $r< 16$ vertices with $r$ even.
    Whenever (A) plays on one of these paths, (I) answers in the same path with the optimal strategy. Since all paths are of even order, the alternating of the moves should not change, and the games on each path results in a (locally) optimal game for (I).
    Therefore, the final score on each $P_{16}$ is at least~$3$, and if $r> 0$ the score on $P_r$ is at least~$1$.
    However, the edges linking any two incident paths may be monochromatic, resulting in a additional score of $-k$ (or $-k+1$ if $r=0$).
    In both cases, we get a total score of $3k - k + 1 = 2k+1 =  2\lfloor\frac{n}{16}\rfloor + 1$.

    The proof for (b) is similar.
    Admirable splits the path into $k$ paths of order~$16$,
    and a remaining path of odd length $r < 16$.
    (A) starts playing an optimal move on $P_r$, then follows each move of (I) with an optimal strategy (either the started strategy on $P_r$,
    or the (I)-start strategy on the $P_{16}$).
    Since the parity of uncolored vertices on each path remains unchanged after each move of
    Admirable, the players alternation is consistent.
    At the end of the game, each path on~$16$ vertices gets a score of at most~$3$,
    and the $P_r$ gets a score of at most~$4$. There are $k$ extra edges linking different paths, that contribute by at most $+k$ to the score. This strategy therefore ensures a score of at most $4k+4 = 4\lceil\frac{n}{16}\rceil$, as desired.
\end{proof}

As a corollary of this lemma and Lemma~\ref{l:imagination} applied on the last vertex $v$ on the path (and so, $\Delta = 2$ and $\deg(v)=1$), we get:

\begin{theorem}\label{thm:paths-A-start}
The (A)-start game balance number of a path $P_n$ of order~$n$ satisfies
\[
\begin{array}{rcccll}
2\left\lfloor \frac{n}{16}\right\rfloor + 1 & \le & \bA(P_n) & \le & 4 \left\lceil \frac{n}{16} \right\rceil  + 5, & \mbox{if $n$ is even;} \1 \\
2\left\lfloor \frac{n}{16}\right\rfloor -4 & \le & \bA(P_n) & \le & 4 \left\lceil \frac{n}{16} \right\rceil, & \mbox{if $n$ is odd.}
\end{array}
\]
\end{theorem}

Note that this in particular answers Question~1 from \cite{original}.
Indeed, even though they are playing the cordiality game, Impish may decide to adopt
our balance game strategy and ensure a score growing linearly with the number of vertices.
As a consequence, as soon as $n\ge 48$, $c_g(P_n) \ge 2$ and the values~$0$ and~$1$ are not assumed anymore by $c_g(P_n)$, answering by the negative to the open question.

In a similar way we obtain:

\begin{lemma} The following holds:
\begin{enumerate}[label=\rm(\alph*)]
    \item For $n$ odd, $\bI(P_n) \ge 2\lfloor\frac{n}{16}\rfloor$.
    \item For $n$ even,  $\bI(P_n) \le 4\lceil\frac{n}{16}\rceil - 1$.
    \end{enumerate}
\end{lemma}

\begin{theorem}\label{thm:paths-I-start}
The (I)-start game balance number of a path $P_n$ of order~$n$ satisfies
\[
\begin{array}{rcccll}
2\left\lfloor \frac{n}{16}\right\rfloor - 5 & \le & \bI(P_n) & \le & 4 \left\lceil \frac{n}{16} \right\rceil - 1, & \mbox{if $n$ is even;} \1 \\
2\left\lfloor \frac{n}{16}\right\rfloor & \le & \bI(P_n) & \le & 4 \left\lceil \frac{n}{16} \right\rceil + 4, & \mbox{if $n$ is odd.}
\end{array}
\]
\end{theorem}

A common substantial point in the proofs of Lemma \ref{l:rec4}, Remark \ref{r:rec6}, and Lemma \ref{l:rec16} is that $\bA(P_n)=\bI(P_n)$ holds for $n=4,6,16$.
Via this approach as applied for $P_{16}$, the following can be proved as a strengthening of the preceding results.
The explicit functions linear in $n$ are obtained by inspecting the game balance numbers of $P_n$, for $n\leq 15$.

\begin{theorem}
For every $n$ we have
$2\leq \bA(P_{n+16})-\bA(P_n)\leq 4$
and
$2\leq \bI(P_{n+16})-\bI(P_n)\leq 4$.
As a consequence, for $n$ even,
\[
 \left\lfloor \frac{n}{8} \right\rfloor + 1 \ \leq \ \bA(P_n) \,, \ \bI(P_n) \ \leq \ 2\left\lceil \frac{n}{8} \right\rceil - 1 ,
\]
while for $n$ odd,
\[
 \left\lfloor \frac{n}{8} \right\rfloor \ \leq \ \bA(P_n) \,, \ \bI(P_n) \ \leq \ 2\left\lceil \frac{n}{8} \right\rceil .
\]
\end{theorem}

Although we do not have a formula with tight asymptotics, we can prove at least that both $\bA(P_n)$ and $\bI(P_n)$ are of the form $cn+o(n)$ for some constant $c$ as $n$ gets large.

\begin{theorem}
The sequences $\bA(P_n)/n$ and $\bI(P_n)/n$ converge to the same value.
\end{theorem}
\begin{proof}
    By \Cref{l:imagination}, it is enough to show that the second sequence converges. To see that this is the case, consider two even integers $m\le n$. The path $P_n$ can be divided into $\lfloor n/m \rfloor$ segments of length $m$ and one additional segment of length $n$ mod $m$ and the player who moves second may play an optimal strategy on each segment of length $m$. It follows that $\bA(P_n)\ge n(\bA(P_m)-1)/m-m$ and $\bI(P_n)\le n(\bI(P_m)+1)/m+m$. However, by \Cref{prop:start-continuity}, $\bI(P_m)\le \bA(P_m) +8$ and $\bA(P_n)\le \bI(P_n)$ so that
    \begin{equation*}
        \frac{\bA(P_m)}{m} - \frac{1}{m} - \frac{m}{n}\le \frac{\bA(P_n)}{n} \le \frac{\bI(P_n)}{n}  \le \frac{\bA(P_m)}{m} + \frac{9}{m} + \frac{m}{n}.
    \end{equation*}
    Thus, it follows that $(\bA(P_{2n})/2n)_{n\in \N}$ and $(\bI(P_{2n})/2n)_{n\in \N}$ converge to the same value. The convergence of the original sequences $(\bA(P_{n})/n)_{n\in \N}$ and $(\bI(P_{n})/n)_{n\in \N}$ follows directly from \Cref{l:imagination}.
\end{proof}

\begin{problem}
    Determine the value of  $\displaystyle{\lim_{n\rightarrow\infty}\bA(P_n)/n = \lim_{n\rightarrow\infty}\bI(P_n)/n}$. 
\end{problem}

\section{Remarks about the cordiality game}

Since the score in the cordiality game is simply the modulus of the score in the balancing game, it follows that $\cA(G)\ge \bA(G)$ and $\cI(G)\ge \bI(G)$ for all graphs $G$. The fact that graphs of odd order and negative balancing number exists (for example, $\bI(C_5)=-1$) shows that these inequalities may be strict, and it is not hard to find further small examples. However, we do not know of a sequence of graphs $(G_k)_k$ such that $\cA(G_k)-\bA(G_k)$ or $\cI(G_k)-\bI(G_k)$ diverge. 

\begin{problem}
    Does there exist a universal constant $C$ such that $\cA(G)-\bA(G)\le C$ and $\cI(G)-\bI(G)\le C$ for all graphs $G$?
\end{problem}

The effectiveness of the greedy strategies presented in the proof of \Cref{thm:global-bounds} lend credibility to an affirmative answer to this question, and in fact, a greedy strategy yields similar global bounds as in \Cref{thm:global-bounds} for the cordiality game.

\begin{theorem}
If $G$ is a graph of order~$n$ with maximum degree $\Delta$, then we have
        \[
    \begin{cases}
            \cA(G) \le \frac{n}{2}+ 2\Delta, &\text{if $n$ is even;}\\
            \cA(G) \le \frac{n}{2} + 3\Delta, &\text{if $n$ is odd;}
        \end{cases}
    \]
    and
            \[
    \begin{cases}
            \cI(G) \le \frac{n}{2} + 3\Delta,  &\text{if $n$ is even;}\\
            \cI(G) \le \frac{n}{2} + 2\Delta, &\text{if $n$ is odd.}
        \end{cases}
    \]
\end{theorem}
\begin{proof}[Proof sketch.]
    We will only consider the (A)-start game for even $n$, as the other cases are similar. To prove the claim, it is enough to show that if $2k$ moves have been played, and the game is in a state $(S_0,S_1)$ such that $\vert s(S_0,S_1)\vert\leq k+2\Delta$, then there is a vertex $v\notin S_0\cup S_1$ for (A) to play such that for all $w\notin S_0\cup S_1 \cup \{v\}$, we have $\vert s(S_0\cup\{v\},S_1\cup\{w\})\vert\leq k+1+2\Delta$.

    Let $(S_0,S_1)$ be such a state, and suppose first that $s(S_0,S_1)$ is non-negative. We choose $v$ such that $s(v;S_0,S_1)$ as defined in the proof of \Cref{thm:global-bounds} is maximal, i.e., such that $s(S_0\cup \{v\},S_1)$ is minimal. As in the proof of \Cref{thm:global-bounds}, we observe that $s(w;S_0\cup \{v\},S_1)\leq s(w;S_0,S_1)+1$, meaning that for all $w$, we have $s(S_0\cup \{v\},S_1\cup \{w\})\leq s(S_0,S_1)+1$. On the other hand, each move by either player can change the score by at most $\Delta$ so that because $s(S_0,S_1)\geq 0$, we have $s(S_0,S_1)\geq -2\Delta$.

    If $s(S_0,S_1)$ is negative, we choose $v$ such that $s(v;S_0,S_1)$ is minimal, i.e., such that $s(S_0\cup \{v\},S_1)$ is maximal. Since $s(w;S_0\cup \{v\},S_1)\geq s(w;S_0,S_1)$ for all $w$, and the score $s(S_0,S_1)<0 -1$, we have that for all $w\notin S_0\cup S_1 \cup \{v\}$,
    \begin{equation*}
        -k-2\Delta\leq s(S_0,S_1) \leq s(S_0\cup \{v\},S_1\cup \{w\})< 2\Delta
    \end{equation*}
    which completes the proof.
\end{proof}

For paths, we can further corroborate the notion that the balance and cordiality game numbers should behave similarly. As for the balance game, we computed $\cA(P_n)$ and $\cI(P_n)$ for $n$ up to 18 (see \Cref{t:cordiality-scores}). While there are differences between the cordiality game numbers and their balance analogues, the values coincide for $n\geq 12$, suggesting that the best strategy in the cordiality game for Impish is to adopt the same strategy as in the balance game. 

\begin{table}[htb]\begin{center}
\begin{tabular}{|l|cccccccccccccccccc|}
 \hline
 $n$ & 1 & 2 & 3 & 4 & 5 & 6 & 7 & 8 & 9 & 10 & 11 & 12 & 13 & 14 & 15 & 16 & 17 & 18 \\
 \hline
 (A)-start & 0 & 1 & 0 & 1 & 2 & 1 & 2 & 3 & 2 & 1 & 4 & 1 & 4 & 1 & 4 & 3 & 4 & 3 \\
 (I)-start & 0 & 1 & 0 & 1 & 2 & 1 & 2 & 3 & 2 & 3 & 2 & 3 & 2 & 3 & 2 & 3 & 2 & 5 \\
\hline
\end{tabular}
\caption{Scores of the cordiality game for short paths}\label{t:cordiality-scores}
\end{center}
\end{table}

Just as in the cordiality game, we can subdivide longer paths into segments of length 16 to find strong strategies for both players. Analogously to \Cref{thm:paths-A-start} and \Cref{thm:paths-I-start}, we obtain the following theorem, and we therefore omit a proof of Theorem~\ref{thm:path-cordiality-game}. We remark that since $\cA\geq \bA$ and $\cI\geq \bI$ in general, the lower bounds in the theorem can be obtained as direct corollaries of \Cref{thm:paths-A-start} and \Cref{thm:paths-I-start}. 

\begin{theorem}
\label{thm:path-cordiality-game}
 For all $n\in \mathbb{N}$, we have
\[
\begin{array}{rcccll}
2\left\lfloor \frac{n}{16}\right\rfloor + 1 & \le & \cA(P_n) & \le & 4 \left\lceil \frac{n}{16} \right\rceil  + 5, & \mbox{if $n$ is even;} \1 \\
2\left\lfloor \frac{n}{16}\right\rfloor -4 & \le & \cA(P_n) & \le & 4 \left\lceil \frac{n}{16} \right\rceil, & \mbox{if $n$ is odd,}
\end{array}
\]
and
\[
\begin{array}{rcccll}
2\left\lfloor \frac{n}{16}\right\rfloor - 5 & \le & \cI(P_n) & \le & 4 \left\lceil \frac{n}{16} \right\rceil - 1, & \mbox{if $n$ is even;} \1 \\
2\left\lfloor \frac{n}{16}\right\rfloor & \le & \cI(P_n) & \le & 4 \left\lceil \frac{n}{16} \right\rceil + 4, & \mbox{if $n$ is odd.}
\end{array}
\]
\end{theorem}

We remark that the lower bound in the theorem above answers Question 1 in \cite{original}, while the upper bound improves Theorem 2.2 in the same article.

\section*{Acknowledgments} 

The authors thank the Faculty of Natural Sciences and Mathematics of the University of Maribor, Slovenia, for hosting the Workshop on Games on Graphs in June 2023.

Research of the second author was supported in part by the South African National Research Foundation (grants 132588, 129265) and the University of Johannesburg. Research of the third author was supported in part by the National
Research, Development and Innovation Office, NKFIH Grant FK 132060. Research of the fourth author was supported by the London Mathematical Society and the Heilbronn Institute for Mathematical Research through an Early Career Fellowship.


\end{document}